\documentclass[10pt,leqno]{article}
\usepackage{amssymb,amsmath,amscd,dsfont}
\usepackage{ifthen}
\newcommand{\tx}[1]{\ensuremath{\text{\upshape #1}}}
\newcommand{\wt}[1]{\ensuremath{\widetilde{#1}}}

\newcommand{\se}{\ensuremath{\subseteq}}
\newcommand{\cP}{\ensuremath{\mathcal{P}}}
\newcommand{\cH}{\ensuremath{\mathcal{H}}}
\newcommand{\Q}{\ensuremath{\mathds{Q}}}
\newcommand{\C}{\ensuremath{\mathds{C}}}
\newcommand{\B}{\ensuremath{\mathds{B}}}
\newcommand{\R}{\ensuremath{\mathds{R}}}
\newcommand{\N}{\ensuremath{\mathds{N}}}
\newcommand{\Z}{\ensuremath{\mathds{Z}}}

\newcommand{\bS}{\ensuremath{\mathds{S}}}
\newcommand{\id}{\ensuremath{\mathop{\text{id}}}}
\newcommand{\SU}{\ensuremath{\tx{SU}_{2}(\C)}}
\newcommand{\Gau}{\ensuremath{\tx{Gau}}}
\newcommand{\Aut}{\ensuremath{\tx{Aut}}}
\newcommand{\Bun}{\ensuremath{\tx{Bun}}}
\newcommand{\Conn}{\ensuremath{\tx{Conn}}}
\newcommand{\pfb}[1][M]{\ensuremath{(K,\pi,P,#1)}}
\newcommand{\from}{\ensuremath{\nobreak:\nobreak}}
\renewcommand{\to}{\ensuremath{\nobreak\rightarrow\nobreak}}
\usepackage[a4paper,twoside,hscale=0.7,vscale=0.82,hmarginratio=3:2,includehead,headheight=14pt]{geometry}
\usepackage{fancyhdr}
\fancypagestyle{plain}{\fancyhead{}\fancyfoot{}}

\newcommand{\shortTitle}{}

\usepackage{amssymb,amsmath,amscd,amsthm,dsfont} 
\usepackage{graphicx}
\usepackage{hyperref}
\usepackage{calc}

\theoremstyle{definition}
\newtheorem{definition}{Definition}[section]
\newtheorem{remark}[definition]{Remark}

\theoremstyle{plain}
\newtheorem{lemma}[definition]{Lemma}
\newtheorem{proposition}[definition]{Proposition}
\newtheorem{theorem}[definition]{Theorem}
\newtheorem{corollary}[definition]{Corollary}
\newtheorem*{nntheorem}{Theorem}
\newenvironment{prf}{\begin{proof}[\textbf{\upshape Proof.}]}{\end{proof}}

\renewcommand{\pfb}[1][M]{\ensuremath{(K,\eta:P\to #1)}}
\newlength\picturewidth
\begin{document}
\sloppy \title{\textbf{The Samelson Product and Rational Homotopy for
Gauge Groups}\\[12pt] \large - preprint -}
\author{Christoph Wockel\\
        Fachbereich Mathematik\\
        Technische Universit\"at Darmstadt\\
\small \texttt{wockel@mathematik.tu-darmstadt.de}}
\renewcommand{\shortTitle}{The Samelson Product and Rational Homotopy}
\maketitle
\thispagestyle{empty}
\begin{abstract}
\noindent This paper is on the connecting homomorphism in the long
exact homotopy sequence of the evaluation fibration
$\tx{ev}_{p_{0}}:C(P,K)^{K}\to K$, where $C(P,K)^{K}\cong\Gau(\cP)$ is
the gauge group of a continuous principal $K$-bundle $P$ over a closed
orientable surface or a sphere.  We show that in this cases the
connecting homomorphism in the corresponding long exact homotopy
sequence is given in terms of the Samelson product. As applications,
we exploit this correspondence to get an explicit formula for
$\pi_{2}(\Gau(\cP_{k}))$, where $\cP_{k}$ denotes the principal
\mbox{$\bS^{3}$-bundle} over $\bS^{4}$ of Chern number $k$ and derive
explicit formulae for the rational homotopy groups
$\pi_{n}(\Gau(\cP))\otimes \Q$.\\[\baselineskip] \textbf{Keywords:}
bundles over spheres, bundles over surfaces, gauge groups, pointed
gauge groups, homotopy groups of gauge groups, rational homotopy
groups of gauge groups, evaluation fibration, connecting homomorphism,
Samelson product, Whitehead product\\[\baselineskip]
\textbf{MSC:} 57T20, 57S05, 81R10, 55P62\\
\textbf{PACS:} 02.20.Tw, 02.40.Re
\end{abstract}
\section*{Introduction}
The topological properties of gauge groups play an important role in
the analysis of the configuration space in quantum field theory.
There one analyses the moduli space $\Conn(\cP)/\Gau(\cP)$ of
connections on a principal $K$-bundle $\cP$ modulo $\Gau(\cP)$ the
group of gauge transformations, shortly called gauge group
(cf. \cite{singer78}). Since $\Conn(\cP)$ is an affine space, the
exact homotopy sequence gives detailed information on the homotopy
groups of the configuration space in terms of the homotopy groups of
the gauge group.  On the other hand, $\pi_{1}(\Gau(\cP))$ and
$\pi_{2}(\Gau(\cP))$ carry crucial information on central extensions
of $\Gau(\cP)$ (cf. \cite{neeb03}), which are important for an
understanding of the relation between the projective and unitary
representations of $\Gau(\cP)$. Furthermore, if $\cP$ is a bundle over
$\bS^{1}$, then $\Gau(\cP)$ is isomorphic to a twisted loop goup, and
thus gauge groups are closely related to Kac-Moody groups
(cf. \cite{mickelsson87}).

We now describe our results in some detail. In the first section, we
recall some basic facts from elementary topology and from the
classification of principal $K$-bundles over spheres and surfaces. The
latter are the types of bundles this text deals with since they have
explicit descriptions in terms of $\pi_{m}(K)$. In the case of a
principal $K$-bundle $\cP=\pfb[\bS^{m}]$ over $\bS^{m}$, this leads to
an explicit description of the gauge group $\Gau(\cP)$ as a subgroup
of $C(\B^{m},K)$ and of $\Gau_{*}(\cP)$ as $C_{*}(\bS^{m},K)$, where
$\Gau_{*}(\cP)$ denotes the group of gauge transformations fixing
$\eta^{-1}(x_{0})$ pointwise and
\mbox{$\B^{m}:=\{x\in\R^{n}:\|x\|_{\infty}\leq 1\}$}.  This
description of $\Gau(\cP)$ leads directly to the main result of this
paper.
\begin{nntheorem}
If $\cP=\pfb[\bS^{m}]$ is a continuous principal $K$-bundle over
$\bS^{m}$, $K$ is locally contractible and $b\in \pi_{m-1}(K)$ is
characteristic for $\cP$, then the connecting homomorphism
\mbox{$\delta_{n}:\pi_{n}(K)\to \pi_{n+m-1}(K)$} in the exact sequence
\[\begin{CD}
\dots
\to \pi_{n+1}(K)
@>{\delta_{n+1}}>>\pi_{n+m}(K)
\to \pi_{n}(\Gau(\cP))
\to\pi_{n}(K)
@>{\delta_{n}}>>\pi_{n+m-1}(K)
\to \dots 
\end{CD}\]
is given by
$\delta_{n}(a)=-\langle a,b\rangle$, where $\langle \cdot ,\cdot
\rangle$ denotes the Samelson product.
\end{nntheorem}
The connection between the Samelson Product and the evaluation
fibration is not new (cf. \cite[Th. 3.2]{whiteheadG.W.46} and
\cite[Sect. 1]{barrattJamesStein60} and Remark
\ref{rem:alternativeProoOfMainTheorem}). The remarkable thing in this
paper is that the above theorem can be proven by using only very
elementary facts on fibrations. However, we give an alternative proof
of the theorem in terms of more involved facts from homotopy theory.

As an application of the above theorem we obtain a new proof of
\cite{kono91} providing an explicit formula for
$\pi_{2}(\Gau(\cP_{k}))$, where $\cP_{k}$ denotes the principal
$\SU$-bundle over $\bS^{4}$ of Chern number $k$. Furthermore, we show
that the connecting homomorphism of the evaluation fibration for
bundles over closed compact orientable surfaces is also given in terms
of the Samelson product, since the situation there reduces to the
situation of bundles over $\bS^{2}$. 

Since the rational Samelson produce \mbox{$\langle \cdot,\cdot
\rangle\otimes \tx{id}_{\Q}$} between the rational homotopy groups
\mbox{$\pi_{n}(K)\otimes \Q$} and \mbox{$\pi_{m}(K)\otimes \Q$}
vanishes for a connected Lie group $K$, this leads to the following
explicit description of the rational homotopy groups of $\Gau(\cP)$ for
a large class of bundles.

\begin{nntheorem}
Let $K$ be a connected Lie group and $\cP =\pfb$ be a continuous
principal $K$-bundle over $\bS^{m}$ or a compact orientable
surface $\Sigma$.
\begin{alignat*}{2}
i)\;&\text{If }M=\bS^{m}&\text{, then }&\pi_{n}^{\Q}(\Gau(\cP))\cong
\pi_{n+m}^{\Q}(K)\oplus \pi_{n}^{\Q}(K).\\
ii)\;&\text{If }M=\Sigma&\text{, then }&\pi_{n}^{\Q}(\Gau(\cP))\cong
\pi_{n+2}^{\Q}(K)\oplus \pi_{n+1}^{\Q}(K)^{2g}\oplus \pi_{n}^{\Q}(K).
\end{alignat*}
\end{nntheorem}
\section*{Acknowledgements} 
The work on this paper was financially supported by a doctoral
scholarship from the Technische Universit\"at Darmstadt.  The author
would like thank Karl-Hermann Neeb and Linus Kramer for giving crucial
hints on the Samelson product. He would also like to thank Mamoru
Mimura and Kouzou Tsukiyama for a very friendly and fruitful
communication out of which grew Remark
\ref{rem:alternativeProoOfMainTheorem}.
\section{General Remarks and Notation}
\begin{remark}
Throughout this paper, we denote by $\B^{n}:=\{x\in
\R^{n}:\|x\|_{\infty}\leq 1\}$ the closed unit ball of radius 1, where
$\|\cdot\|_{\infty}=\max\{|x_1|,\ldots,|x_{n}|\}$ denotes the
infinity-norm (we use this somewhat uncommon setting since then the
proof of Theorem \ref{thm:connectingHomomorphismIsTheSamelsonProduct}
becomes less cryptic).  Furthermore we set $I=[-1,1]=\B^{1}$ and thus
have $\B^{n}=\B^{n-1}\times I$. By $\bS^{n}$, we denote the $n$-sphere
and identify it interchangeably with \mbox{$\{x\in
\R^{n+1}:\|x\|=1\}$} (where $\|\cdot\|$ denots the euclidean norm),
with \mbox{$\{x\in \R^{n+1}:\|x\|_{\infty}=1\}$} or with
\mbox{$\B^{n}/\partial \B^{n}$}, depending on what is convenient in
the considered situation. When dealing with pointed spaces, we take
$(1,0,\ldots,0)$ as the base-point in $\B^{n}$, \mbox{$\{x\in
\R^{n+1}:\|x\|=1\}$} or \mbox{$\{x\in \R^{n+1}:\|x\|_{\infty}=1\}$}
and $\partial \B^{n}$ as base-point in $\B^{n}/\partial \B^{n}$.

If $\sim$ is an equivalence relation on the topological space $X$ and
$X/\sim$ is the quotient $X$ by this relation, then the continuous
functions on $X/\sim$ are in on-to-one correspondence with the
continuous functions on $X$ which are constant on the equivalence
classes of $\sim$ \cite[\S I.3.4]{bourbakiTop}.

If $f\from X\times Y\to Z$ is a function, then we denote for each $x\in
X$ by $f_{x}$ the function $f_{x}\from Y\to Z$, $y\mapsto
f(x,y)$, and for each $y\in Y$ by $f_{y}$ the function $f_{y}\from X\to
Z$, $y\mapsto f(x,y)$.

If $X$,$Y$ are spaces with base-points $x_{0}$, $y_{0}$,
then \mbox{$C_{*}(X,Y):=\{f\in
C(X,Y):f(x_{0})=y_{0}\}$}. If $X=I$ we set
$PY:=C_{*}(I,Y)$ and if $X=\bS^{1}$ we set $\Omega
Y:=C_{*}(\bS^{1},Y)$
\end{remark}
\begin{remark}
If $X$ and $Y$ are topological spaces, then we equip $C(X,Y)$ with the
compact-open topology. If $Y=K$ is a topological group, then the
compact-open topology on $C(X,K)$ coincides with the topology of
compact convergence (cf. \cite[Th. X.3.4.2]{bourbakiTop}) and this
turns $C(X,K)$ into a topological group.

The elementary facts on the compact open topology on $C(X,K)$ we use
throughout this paper are the following (cf. \cite{bourbakiTop}):
\begin{itemize}
\item If $x\in X$, then the \textit{exaluation map}
        $\tx{ev}_{x}:C(X,Y)\to Y$ is continuous.
\item If $Z$ is a topological space and $f\from X\to Z$ is continuous, then
	the \textit{pull-back} $f^{*}\from C(Z,Y)\to C(X,Y)$,
	$\gamma \mapsto \gamma \circ f$ is continuous.
\item We have a continuous map $~^{\wedge}\from C(X\times Y,Z)\to C(X,C(Y,Z))$,
	$f^{\wedge}(x)(y)=f(x,y)$ and if $Y$ is locally compact, then
	this map is a homeomorphism. This property is called the
	\textit{exponential law} or \textit{Cartesian closedness principle}.
\end{itemize}
\end{remark}
\begin{remark}\label{rem:classificationOfBundlesOverSpheres}
If $m\in \N^{+}$, then the equivalence classes of continuous principal
$K$-bundles over $\bS^{m}$ are in one-to-one correspondence with the
orbits of the $\pi_{0}(K)$-action on $\pi_{m-1}(K)$, where
$\pi_{0}(K)$ acts on $\pi_{m-1}(K)$ by $([\gamma],k)\mapsto [k\gamma
k^{-1}]$ (cf. \cite[\S 18.5]{steenrod51}).

A characteristic map for a fixed bundle $\cP=\pfb[\bS^{m}]$ can be
obtained as follows. Take $\bS^{n}:=\{x\in \R^{n}:\|x\|=1\}$ and
$\bS^{n-1}=\{x\in \R^{n-1}:\|x\|=1\}=\bS^{n}\cap
\{x\in\R^{n}:x_{n}=0\}$ and let $V_{N/S}\se\bS^{n}$ denote open
$n$-cells with $\bS^{n-1}\se V_{N/S}$ and $(0,\ldots,0,1)\in V_{N}$
and $(0,\ldots,0,-1)\in V_{S}$. Then there exist sections
$\sigma_{N/S}:V_{N/S}\to P$ and $\sigma_{S}(x)=\sigma_{N}(x)\gamma
(x)$ defines a continuous map $\gamma :\bS^{n-1}\to K$. If we
substitute $\sigma_{N}$ by $\sigma_{N}\cdot \gamma (x_{0})$ we may
assume that $\gamma (x_{0})=e$. Then $[\gamma]\in \pi_{n-1}(K)$
represents the equivalence class of $\cP$ (cf. \cite[\S 18.1]{steenrod51}).
\end{remark}
\begin{remark}\label{rem:classificationOfBundlesOverSurfaces}
Let $K$ be a connected topological group and $\Sigma$ be a closed
compact orientable surface. For the set of equivalence classes of
continuous principal $K$-bundles over $\Sigma $ we have that it is
equal to $[\Sigma,BK]$, where $BK$ is the classifying space of $K$
(cf. \cite[Th. 4.13.1]{husemoller66}). Furthermore we have
\[
[\Sigma,BK]\cong H^{2}(\Sigma,\pi_{2}(BK))\cong
\tx{Hom}(H_{2}(\Sigma),\pi_{1}(K))\cong \pi_{1}(K).
\]
The first isomorphism is a
consequence of \cite[Cor. VII.13.16]{bredon} and \cite[Th.
VII6.7]{bredon}, the second is \cite[Th. V.7.2]{bredon} which
applies since $H_{1}(\Sigma)\cong\Z^{2g}$ is free, and the last
isomorphism follows from $H_{2}(\Sigma)\cong \Z$.
\end{remark}
\begin{remark}\label{rem:constructionOfTheConnectingHomomorphism}
We recall the construction of the connecting homomorphism for a
fibration $p\from Y\to B$ with fibre $F=p^{-1}(\{x_{0}\})$.
This fibration yields a long exact homotopy sequence
\[\begin{CD}
\ldots
\to
\pi_{n+1}(B)
@>{\delta_{n+1}}>>
\pi_{n}(F)
@>{\pi_{n}(i)}>>
\pi_{n}(Y) 
@>{\pi_{n}(q)}>>
\pi_{n}(B)
@>{\delta_{n}}>>
\pi_{n-1}(F)\to \ldots
\end{CD}\]
and the construction of the connecting homomorphism $\delta_{n}$ is as
follows (cf. \cite[Th. VII.6.7]{bredon}): Let ${f}\in
C_{*}(\B^{n},B)$ represent an element of $\pi_{n}(B)$,
i.e. $\left.f\right|_{\partial \B^{n}}\equiv x_{0}$. Then
$f$ can be lifted to a map
$F:\B^{n}\to Y$ with $q\circ F=f$ since $q$ is a fibration. Then
$F$ takes $\partial \B^{n}\cong \bS^{n-1}$ into $q^{-1}(x_{0})=F$,
and $\left.F\right|_{\partial \B^{n}}$ represents $\delta ([{f}])$.
\end{remark}
\section{The Connecting Homomorphism}
\begin{definition}[Bundle Map, Automorphism Group, Gauge Group]
If $\cP =\pfb$ and $\cP '=(K,\eta '\from P'\to B')$ are principal $K$-bundles,
then 
\[
\Bun(\cP ,\cP'):=\{f\in C(P,P'):(\forall p\in P)(\forall k\in K)\;
f(p\cdot k)=f(p)\cdot k\}
\]
are called \textit{bundle maps} from $\cP$ to $\cP'$. Furthermore,
\mbox{$\Aut(\cP):=\Bun(\cP,\cP)\cap \tx{Homeo}(P)$} is called the
group of bundle automorphism or \textit{automorphism group} of $\cP$
and \mbox{$\Gau(\cP):=\{f\in \Aut(P):\eta \circ f=\eta\}$} is called
the group of bundle equivalences or \textit{gauge group} of $\cP$.
\end{definition}
\begin{remark}
The gauge group of $\cP$ is isomorphic to the space of continuous 
$K$-equivariant mappings 
\[
C(P,K)^{K}:=\{f\in C(P,K):(\forall p\in P)(\forall k\in K)\;
f(p\cdot k)=k^{-1}\cdot f(p)\cdot k\}
\]
under the isomorphism $C(P,K)^{K}\ni f\mapsto \big(p\mapsto p\cdot
f(p)\big)\in \Gau(\cP)$, and we endow $C(P,K)^{K}$ with the subspace
topology induced from the compact-open topology on $C(P,K)$. This
turns $C(P,K)^{K}$ and thus $\Gau(\cP)$ into topological groups.
\end{remark}
\begin{remark}\label{rem:descriptionOfBundlesOverSpheres}
We recall the description of principal $K$-bundles over $\bS^{m}$ by
its characteristic maps (also called clutching functions). Given a
principal $K$-bundle $\cP=\pfb[\bS^{m}]$ over $\bS^{m}$ and denoting
by $q:\B^{m}\to \bS^{m}$ the quotient map identifying $\partial
\B^{m}$ with the base-point in $\bS^{m}$, \cite[Cor.
VII.6.12]{bredon} provides a map $\sigma :\B^{m}\to P$ satisfying
$\eta \circ \sigma=q$.  Thus $\sigma (y)\cdot \gamma (y)=\sigma
(y_{0})$ for $y\in \partial \B^{m}$ and we thus obtain a continuous
map $\gamma :\partial \B^{m}\cong \bS^{m}\to K$ satisfying $\gamma
(y_{0})=e$ which is called the \textit{clutching function} or
\textit{characteristic map} describing $\cP$.  Furthermore $\gamma$ is
a representative of $[\cP]$ since we may identify $\tx{int}(\B^{m})$
with $V_{N}$, $(\B^{m}/\partial\B^{m})\backslash \{0\}$ with $V_{S}$
and then
\begin{alignat*}{2}
 \sigma_{N}  &\from V_{N}\to P,\;\; & x\mapsto \sigma (x)&\\
 \sigma_{S}  &\from V_{S}\to P,\;\; & x\mapsto \sigma (x)\cdot&
   \gamma \left(\frac{x}{\|x\|_{\infty}}\right)
\end{alignat*}
denote corresponding sections (cf. Remark
\ref{rem:classificationOfBundlesOverSpheres}).  Set $P/\gamma
:=\B^{m}\times K/\sim$ with $(x,k)\sim(y,k'):\Leftrightarrow x,y\in
\partial B^{m}$ and $\gamma (x)\cdot k = \gamma (y)\cdot k'$ and endow it with
the quotient topology. Then $K$ acts continuously on $P/\gamma$ by
$([(x,k)],k)'\mapsto [(x,kk')]$ and
\[
P/\gamma \to P,\;\; [(x,k)]\mapsto \left\{
\begin{array}{l@{\;\tx{ if }\;}l}
\sigma_{N}(x)\cdot k & x\in V_{N}\\
\sigma_{S}(x)\cdot \gamma (\frac{x}{\|x\|_{\infty}})\cdot k &
x\in V_{S}
\end{array} \right.
\]
is an isomorphism between the $K$-spaces $P$ and $P/\gamma$, whence a bundle
isomorphism.
\end{remark}
\begin{lemma}\label{lem:isoOfTheGaugeGroups}
Let $\cP=\pfb[\bS^{m}]$ be a continuous principal $K$-bundle with
characteristic map $\gamma :\partial \B^{m}\to K$ and set
\[
D^{m}=(I\times \partial \B^{m-1})\cup (\{1\}\times
\B^{m-1})\cup\{(t,x)\in I\times \R^{m-1}:t=-1\tx{ and }
\frac{1}{2}\leq\|x\|_{\infty}\leq 1\}\se \partial \B^{m}
\]
if $m\geq 2$
\begin{figure}
\centering
\includegraphics[width=\picturewidth*\real{0.71},keepaspectratio]{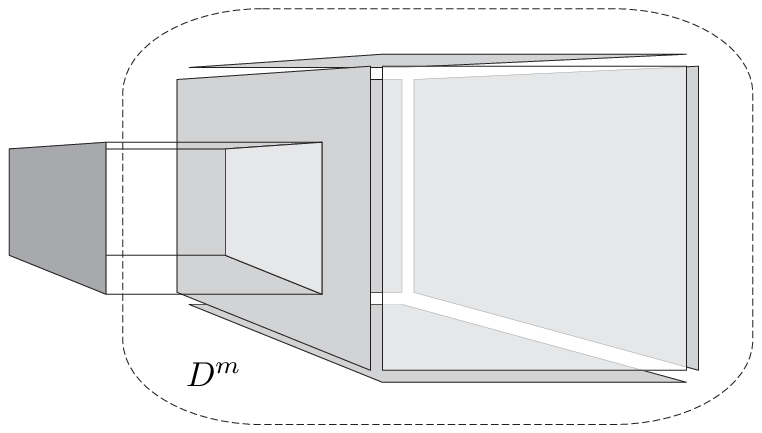}
\caption{Illustration of $D^{m}$}
\label{GP}
\end{figure} (cf. Figure \ref{GP}) and
\mbox{$D^{1}=\{1\}$}. Then
\begin{multline*}
C(P,K)^{K}\cong G(\cP):=\{f\in C(\B^{m},K):
(\exists k\in K)\left.f\right|_{D^{m}}\equiv k, \\
(\forall x\in \partial \B^{m}\backslash D^{m})\;\gamma
(x)^{-1}\cdot f(x)\cdot \gamma (x)=f(x_{0})\}.
\end{multline*}
and thus $G_{*}(\cP):=\{f\in G(\cP):f(x_{0})=e\}\cong C_{*}(\bS^{m},K)$.
\end{lemma}
\begin{prf}
Let $\gamma$ be determined by $\sigma :\B^{m}\to P$ with $\gamma
(x_{0})=e$ as in the preceding remark.  Since $\partial
\B^{m}=(I\times \partial \B^{m-1})\cup( \{-1,1\}\times \B^{m-1})$ and
since $ D^{m}$ is contractible in $\partial \B^{m}$, \cite[Prop.
0.17]{hatcher02} implies that $\gamma $ is homotopic to a map which is
the identity on $D^{m}$. Since homotopic maps yield equivalent bundles
(cf. \cite[Th. 18.3]{steenrod51}) we may assume that
$\left.\gamma\right|_{D^{m}}\equiv e$. Furthermore, $f\mapsto
f\circ\sigma $ provides a map \mbox{$\sigma^{*} \from C(P,K)^{K}\to
G(\cP)$} since $\left.\gamma\right|_{D^{m}}\equiv e$. We claim that
$\sigma^{*}$ is an isomorphism and that an inverse map can be
constructed with $\sigma_{N}$ and $\sigma_{S}$ in terms of pull-backs
and multiplication in spaces of continuous mappings. In fact, for
$f\in G(\cP)$ we set \mbox{$f_{N}:=\left.f\right|_{V_{N}}$} and $f_{S}\from
V_{S}\to K$, $x\mapsto \gamma (\frac{x}{\|x\|_{\infty}})^{-1} f(x)
\gamma (\frac{x}{\|x\|_{\infty}})$. Furthermore,
\mbox{$p=\sigma_{N/S}(\eta (p))\cdot k_{N/S}(p)$} determines
continuous maps $k_{N/S}\from\eta^{-1}(V_{S/N})\to K$ satisfying
\mbox{$k_{N}(p)=\gamma \left(\frac{\eta (p)}{\|\eta
(p)\|_{\infty}}\right)k_{S}(p)$} for \mbox{$p\in \eta^{-1}(V_{S}\cap
V_{N})$}.  Then
\[
f'\from P\to K,\;\;
p\mapsto k_{N/S}(p)^{-1}f_{N/S}(\eta(p))k_{N/S}(p)\;\tx{ if }\;
\eta (p)\in
V_{N/S}
\]
determines an element of $C(P,K)^{K}$ and the assignment $f\mapsto f'$
defines a continuous inverse of $\sigma^{*}$.
\end{prf}
\begin{remark}(cf. \cite[3.7]{pressleysegal86})
Note that the preceding lemma implies that if $\cP=\pfb[\bS^{1}]$ is
a principal $K$-bundle over the circle given by $[k]\in \pi_{0}(K)$,
then the gauge group is isomorphic to the twisted loop group
\[
C_{k}(\bS,K):=\{f\in C(\R,K):f(x+n)=k^{-n}f(x)k^{n}\}.
\]
In fact, since a characteristic map for a bundle over $\bS^{1}$ is
represented by an element $k\in K$ we have \mbox{$G(\cP)=\{f\in
C(I,K):k^{-1}\cdot f(-1)\cdot k= f(1)\}$} and the isomorphism
\[
G(\cP)\ni f\mapsto \left(x\mapsto k^{-1}\cdot f(2(x-n)-1)\cdot k^{n} 
\right)\in C_{k}(\bS,K),
\]
where $n\in\Z$ such that $x-n \in[0,1]$.
\end{remark}
\begin{definition}[Evaluation Map] If $\cP=\pfb[\bS^{m}]$ is a
continuous principal $K$-bundle, then
$\tx{ev}_{x_{0}}:G(\cP)\to K$, $f\mapsto
f(x_{0})$ is called the \textit{the evaluation map}.
\end{definition}
\begin{lemma}
If $\cP=\pfb[\bS^{m}]$ is a continuous principal $K$-bundle and $K$
is locally contractible, then the evaluation map is a fibration
with kernel $G_{*}(\cP)\cong C_{*}(\bS^{m},K)$.
Furthermore, $K_{\cP}:=\tx{im}(\tx{ev}_{x_{0}})$ is open and thus
contains the identity component $K_{0}$.
\end{lemma}
\begin{prf}
Since $K$ is locally contractible, there exist open unit
neighbourhoods $V\se U$ and a continuous map $F:[0,1]\times V\to U$ such
that $F(0,k)= e$, $F(1,k)=k$ for all $k\in V$ and $F(t,e)=e$ for all
$t\in[0,1]$. For $k\in V$ we set $\tau_{k}:=F(\cdot ,k)$, which is a
continuous path and satisfies $\tau_{k}(0)=e$ and $\tau_{k}(1)=k$.
Furthermore, the map $V\ni k\mapsto \tau_{k}\in C(I,K)$ is continuous as an
easy calculation in the topology of compact convergence shows.

Now $V\ni k\mapsto f_{\tau_{k}}\in G(\cP)$ defines a continuous
section of the evaluation map and since $\tx{ev}_{x_{0}}$ is surjective
this shows that $(G_{*}(\cP),\tx{ev}_{x_{0}}:G(\cP)\to K)$ is a
continuous principal $G_{*}(\cP)$-bundle and thus a fibration
(cf. \cite[Cor. VII.6.12]{bredon}). Since the bundle projection of a
locally trivial bundle is open it follows in particular that $\tx{ev}_{x_{0}}$
is open and thus that its image is open.
\end{prf}
\begin{lemma}
If $\cP=\pfb[\bS^{m}]$ is a continuous principal $K$-bundle over
$\bS^{m}$ and $K$ is locally contractible, then the evaluation map
$\tx{ev}_{x_{0}}$ induces a long exact homotopy sequence
\begin{equation}\label{eqn:longExactHomotopySequence}
\begin{CD}\dots
\to \pi_{n+1}(K)
@>{\delta_{n+1}}>>\pi_{n+m}(K)
\to \pi_{n}(\Gau(\cP))
\to 
	\pi_{n}(K)
@>{\delta_{n}}>>\pi_{n+m-1}(K)
\to \dots 
\end{CD}\end{equation}
\end{lemma}
\begin{prf}
Since $K_{\cP}$ contains the identity component $K_{0}$ we have
$\pi_{n}(K_{0})=\pi_{n}(K_{\cP})=\pi_{n}(K)$, and since
\mbox{$\pi_{n+m}(K)=\pi_{0}(C_{*}(\bS^{n+m},K))\cong
\pi_{0}(C_{*}(\bS^{n},C_{*}(\bS^{m},K)))=\pi_{n}(C_{*}(\bS^{m},K))$}
this a direct consequence of the long exact homotopy sequence
(cf. \cite[Th. VII.6.7]{bredon}) for \mbox{$\tx{ev}_{x_{0}}\from
G(\cP)\cong C(P,K)^{K}\cong\Gau(\cP)\to K_{\cP}$} and the
preceding lemma.
\end{prf}
\begin{definition}[Samelson Product]
If $K$ is a topological group, $a\in \pi_{n}(K)$ is represented by
$\alpha \in C_{*}(\bS^{n},K)$ and $b\in \pi_{m}(K)$ is represented by
$\beta \in C_{*}(\bS^{m},K)$, then the commutator map
\[
\alpha\#\gamma:\bS^{n}\times \bS^{m}\to K,\;\;(x,y)\mapsto \alpha
(x)\beta (y)\alpha (x)^{-1}\beta (y)^{-1}
\]
maps $\bS^{n} \vee \bS^{m}$ to $e$. Hence it may be viewed as an element of
$C_{*}(\bS^{n}\wedge\bS^{m},K)$ and thus determines an element
$\langle a,b \rangle:=[\alpha\#\beta]\in
\pi_{0}(C_{*}(\bS^{n+m},K))\cong \pi_{n+m}(K)$. The map
\[
\pi_{n}(K)\times \pi_{m}(K)\to \pi_{n+m}(K),\;\; (a,b)\mapsto \langle
a,b \rangle
\]
is biadditive \cite[Th. X.5.1]{whitehead78} and is called the
\textit{Samelson Product} (cf. \cite[Sect. X.5]{whitehead78}).
\end{definition}
\begin{theorem}\label{thm:connectingHomomorphismIsTheSamelsonProduct}
If $\cP=\pfb[\bS^{m}]$ is a continuous principal $K$-bundle over
$\bS^{m}$, $K$ is locally contractible and $b\in \pi_{m-1}(K)$ is
characteristic for $\cP$ (cf. Remark
\ref{rem:classificationOfBundlesOverSpheres}), then the connecting
homomorphism $\delta_{n}:\pi_{n}(K)\to \pi_{n+m-1}(K)$ in
\eqref{eqn:longExactHomotopySequence} is given by
$\delta_{n}(a)=-\langle a,b\rangle$, where $\langle \cdot ,\cdot
\rangle$ denotes the Samelson product.
\end{theorem}
\begin{prf}
Let $b$ be represented by $\gamma\in C_{*}(\partial
\B^{m},K)$ with $\left.\gamma\right|_{D^{m}}\equiv e$, $a\in
\pi_{n}(K)$ be represented by $\alpha\in C(\B^{n},K)$ with
$\left.\alpha\right|_{\partial \B^{n}}\equiv e$.
Due to the construction of the connecting homomorphism (cf. Remark
\ref{rem:constructionOfTheConnectingHomomorphism}), we have to
construct a lift $A:\B^{n}\to G(\cP)$ of $\alpha$. 

We set $\wt{\alpha}(x,s,t):=\alpha
(x,\frac{t+1}{2}s-(1-\frac{t+1}{2}))$ and note that
$\wt{\alpha}(x,s,1)=\nobreak\alpha (x,s)$ and
$\wt{\alpha}(x,s,-1)=\nobreak\alpha (x,-1)=\nobreak e$.
If $m=1$, then $[\gamma]=[k]\in\pi_{0}(K)$ for some $k\in K$, and we set 
\[
A\from \B^{n}\times I\times I\to K,\;\;(x,s,t)\mapsto 
\wt{\alpha}(x,s,-t)\cdot k\cdot \wt{\alpha }(x,s,t)\cdot k^{-1}.
\]
If $m\geq 2$ the construction of $A$ is as follows. First we set
\[
A':\B^{n-1}\times \B^{2}\times \frac{1}{2}\B^{m-1}\to K,\;\;
(x,s,t,y)\mapsto \wt{\alpha}(x,s,-t)
\wt{\gamma}(y)\wt{\alpha}(x,s,t)\wt{\gamma}(y)^{-1},
\]
where $\wt{\gamma}(y):=\gamma (-1,y)$
\begin{figure}
\centering
\includegraphics[width=\picturewidth*\real{1.16},keepaspectratio]{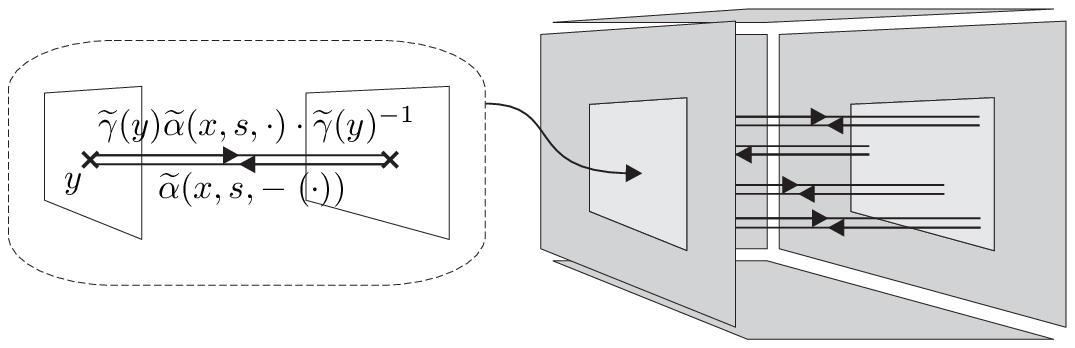}
\caption{Construction of $A'$}
\label{A_prime}
\end{figure} (cf. Figure \ref{A_prime}). Note that due to
$\left.\gamma\right|_{ I\times \partial\B^{m-1}\cup \{1\}\times\B^{m-1}}
\equiv e$, we have that $\wt{\gamma}\from \B^{m-1}\to K$
represents the same element of $\pi_{m-1}(K)$ as $\gamma$ does if we
identify $\bS^{m-1}$ with $\B^{m-1}/\partial \B^{m-1}$ instead of
$\partial \B^{m}$.

Then $t\mapsto A'_{x,s,y}(t)$ satisfies
$A'_{x,s,y}(t)=\nobreak\wt{\alpha}(x,s,-t)\wt{\alpha
}(x,s,t)$ if $\|y\|_{\infty}=\frac{1}{2}$ since then
$\wt{\gamma}(y)=e$ and this map is homotopic to the map which is
constantly $\alpha (x,s)$.  We take a standard homotopy $F_{x,s}$ between
$t\mapsto \wt{\alpha} (x,s,-t)\cdot \wt{\alpha}(x,s,t)$ and the
constant map $\alpha (x,s)$.%

Then $(x,s,r,t)\mapsto F_{x,s}(r,t)$ is continuous and thus
\[
A:\B^{n-1}\times \B^{2}\times \B^{m-1}\to K,\;\;
(x,s,t,y)\mapsto \left\{
\begin{array}{l@{\;\tx{ if }\;}l}
A'(x,s,t,y) & \|y\|_{\infty}\leq \frac{1}{2}\\
F_{x,s}(3-4\|y\|_{\infty},t) & \|y\|_{\infty}\geq\frac{1}{2}
\end{array}
\right.
\]
defines a continuous map
\begin{figure}
\centering
\includegraphics[width=\picturewidth,keepaspectratio]{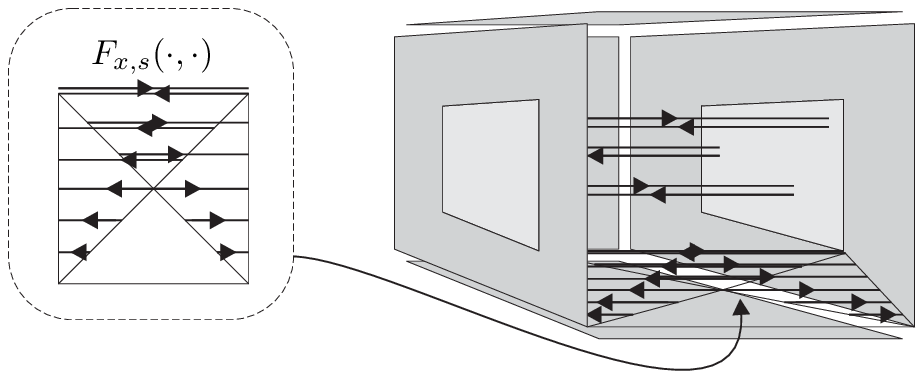}
\caption{Construction of $A$}
\label{A}
\end{figure} (cf. Figure \ref{A}) such that $A_{x,s}$ is an
element of $G(\cP)$ (note that $F_{x,s}$ satisfies
$\left.F_{x,s}\right|_{I\times \{-1,1\}}\equiv \alpha (x,s)$).
Furthermore $(x,s)\mapsto A_{x,s}$ is a lift of $\alpha$ since it is
continuous by the exponential law and satisfies
$A_{x,s}(1,0,\ldots,0)=\alpha (x,s)$.

We now restrict the lift to $\partial \B^{m}=\partial \B^{m-1}\times
I\cup \B^{m}\times \{-1,1\}$. For $x\in \partial \B^{m-1}$ or $s=-1$
we see that $F_{x,s}\equiv e$ since then $\wt{\alpha}(x,s,t)=e$
and thus that in this case $A_{x,s}\equiv e$. Identifying $\bS^{n-1}$
with $\{x\in \partial \B^{n}:x_{n}=1\}$ modulo boundary it thus
suffices to evaluate the lift for $s=1$. Note that we have
$\wt{\alpha}(x,1,t)=\alpha (x,t)$. If $m=1$ we take a homotopy
$G:I\times I\to K$ between $t\mapsto \alpha (x,-t)$ and $t\mapsto
\alpha (x,t)^{-1}$. Then $(r',x,t,y)\mapsto G(r',t)A(x,1,t,y)$ defines
a homotopy in $G_{*}(\cP)$ between $\left.A\right|_{s=1}$ and
$\alpha^{-1}\#\wt{\gamma}$.

If $m\geq 2$, we define $\wt{F}_{x}:I^{3}\to K$, with
$\wt{F}_{x}(1,r,t)=\wt{F}_{x}(r,1,t)=F_{x,1}(r,t)$,
constant on straight lines joining $(1,r,t)$ with $(r,1,t)$ and
$\wt{F}_{x}(r',r,t)=e$ if $r'+r\leq 0$. Then
$\wt{F}_{x}(1,1,t)=F_{x,1}(1,t)=\alpha (x,-t)\alpha (x,t)$,
$\left.F_{x}\right|_{\{-1\}\times I\times I}\equiv e$ and
$F_{x}$ depends continuously on $x$.
Thus
\begin{multline*}
G:I\times \B^{n-1}\times I\times \B^{m-1},\;\;\\(r',x,t,y)\mapsto 
\left\{\begin{array}{l@{\;\tx{ if }\;}l}
\wt{F}_{x}(r',1,t)\alpha(x,t)^{-1}\wt{\gamma}(y)
\alpha (x,t)\wt{\gamma}(y)^{-1}&
      \|y\|_{\infty}\leq \frac{1}{2}\\
\wt{F}_{x}(r',3-4\|y\|_{\infty},t)&
      \|y\|_{\infty}\geq \frac{1}{2}
\end{array} \right.
\end{multline*}
defines a homotopy in $G_{*}(\cP)$ between
$\left.A\right|_{s=1}=G_{1}$ and $\alpha^{-1} \#\wt{\gamma}
=G_{-1}$. Thus we have
\mbox{$[\alpha^{-1}\#\wt{\gamma}]=[\alpha^{-1}\#\gamma]=-\langle a,b
\rangle$} in $\pi_{n+m-1}(K)$.
\end{prf}
\begin{remark}\label{rem:alternativeProoOfMainTheorem}
The above sequence can also be obtained as follows. Let $\cP_{K}
=(K,\eta_{K}\from EK\to BK)$ be a universal bundle for $K$, i.e.
a continuous principal $K$-bundle such that $\pi_{n}(EK)$ vanishes
for $n\in\N^{+}$. Furthermore let \mbox{$\gamma :\bS^{m}\to BK$} be a
classifying map for $\cP$ and denote by $\Gamma \from P\to EK$ the
corresponding bundle map.

Now each $f\in \Bun(\cP,\cP_{K})$ induces a map $\bar{f}:\bS^{m}\to BK$
and the map
\[
\Bun(\cP,\cP_{K},\Gamma )\ni f\mapsto \bar{f}\in C(B,BK,\gamma )
\]
is a fibration \cite[Prop. 3.1]{gottlieb72}, where $\Bun(\cP
,\cP_{K},\Gamma )$ (respectively $C(B,BK,\gamma )$) denotes the
connected component of $\Gamma $ (respectively $\gamma
$). Then the fibre \mbox{$F=\{\Bun(\cP,\cP_{K}):\bar{f}=\gamma \}$} of
this map is homeomorphic to $\Gau(\cP)$ \cite[Prop.
4.3]{gottlieb72}. Since $\Bun(\cP ,\cP_{K})$ is essentially
contractible \cite[Th. 5.2]{gottlieb72}
$\pi_{n}(\Bun(\cP,\cP_{K}))$ vanishes, and thus the long exact
homotopy sequence of the above fibration leads to
\mbox{$\pi_{n-1}(\Gau(\cP))\cong \pi_{n}(C(B,BK,\gamma ))$}
(cf. \cite[Th. 1.5]{tsukiyama85}).

We now consider the evaluation map
\mbox{$\tx{ev}_{x_{0}}:C(\bS^{m},BK)\to BK$} in the base-point $x_{0}$
of $\bS^{m}$. This map is a fibration (\cite[Th. VII.6.13]{bredon})
and we thus get a long exact homotopy sequence
\begin{multline}\label{eqn:longExactSequenceForClassifyingSpace}
\begin{CD}
\ldots\to\pi_{n+1}(BK)@>{\delta_{n+1}}>>\pi_{n}(C_{*}(\bS^{m},BK,\gamma ))\to 
\pi_{n}(C(\bS^{m},BK,\gamma ))\end{CD}\\
\begin{CD}
\to \pi_{n}(BK)@>{\delta_{n}}>>\pi_{n-1}(C_{*}(\bS^{m},BK,\gamma ))\to\ldots
\end{CD}
\end{multline}
If we identify $\pi_{n}(C_{*}(\bS^{m},BK,\gamma))$ with
$\pi_{n+m}(BK)$ (cf.\cite[2.10]{whiteheadG.W.46}), then the connecting
homomorphism in this sequence is given by $\delta_{n+1}(a)=-[a,b]$,
where $b=[\gamma]\in \pi_{m}(BK)$ and $[\cdot,\cdot]$ denotes the
Whitehead product (cf. \cite[Th. 3.2]{whiteheadG.W.46} and
\cite[(3.1)]{whiteheadJ.H.53}).

Since $\pi_{n}(EK)$ vanishes, the connecting homomorphism
$\Delta\from\pi_{n+1}(BK)\to \pi_{n}(K)$ from the long exact homotopy
sequence for $\cP_{K}$ is an isomorphism. Since we have 
\[
\Delta\left([a,b]\right)=(-1)^{n}\langle \Delta (a),\Delta (b) \rangle
\]
for $a\in\pi_{n+1}(BK)$ by \cite[Sect. 1]{barrattJamesStein60},
\eqref{eqn:longExactSequenceForClassifyingSpace} yields a long exact
sequence
\[
\begin{CD}
\ldots\pi_{n}(K)@>{\delta'_{n}}>>\pi_{n+m-1}(K)\to 
\pi_{n-1}(\Gau(\cP))
\to \pi_{n-1}(K)@>{\delta'_{n-1}}>>\pi_{n+m-2}(K)\to\ldots
\end{CD}.
\]
with connecting homomotphism $\delta'_{n}(a)=(-1)^{n}\langle a,b
\rangle$ if we identify $\pi_{n-1}(\Gau(\cP))$ with
$\pi_{n}(C_{*}(\bS^{m},BK,\gamma))$ as described above and
$\pi_{n+1}(BK)$ with $\pi_{n}(K)$ and $\pi_{n+m}(BK)$ with
$\pi_{n+m-1}(K)$ by $\Delta$.
\end{remark}
\section{Applications}
\begin{proposition}(cf. \cite[Lem. 1.3]{kono91})
If $\cP_{k}$ is a principal $\tx{SU}_{2}(\C)$-bundle over $\bS^{4}$ of Chern
number $k\in \Z$, then
$\pi_{2}(\Gau(\cP_{k}))\cong \Z_{\tx{gcd}(k,12)}$. In particular, if
$\cP_{1}=\cH$ is the quaterionic Hopf fibration, then
$\pi_{2}(\tx{Gau}(\cH))$ vanishes.
\end{proposition}
\begin{prf}
Since by \cite[Th. 6.4.2]{naber00} $\cP_{k}$ is classified by its
Chern number $k\in \Z\cong\pi_{3}(\SU)$, Theorem
\ref{thm:connectingHomomorphismIsTheSamelsonProduct} provides an exact
sequence
\[\begin{CD}
\ldots\to\pi_{3}(\SU)
@>{\delta^{k}_{2}}>>
\pi_{6}(\SU)
@>{\pi_{2}(i)}>> \pi_{2}(\Gau(\cP_{k}))
\to
\pi_{2}(\SU) \to\ldots\;\;,
\end{CD}\]
where $\delta^{k}_{2}:\pi_{3}(\SU)\to \pi_{6}(\SU)$ is given by
$a\mapsto -\langle a ,k\rangle$. Since $\pi_{3}(\SU)\cong \Z$,
$\pi_{6}(\SU)\cong \Z_{12}$ and $\langle 1,1 \rangle$ generates
$\Z_{12}$, we may assume that $\delta^{k}_{2} :\Z\to \Z_{12}$ is the map
$\Z\ni z\mapsto -[kz]\in \Z_{12}$ due to the biadditivity of $\langle
\cdot,\cdot \rangle$. Since $\pi_{2}(\SU)$ is trivial we have that
$\pi_{2}(i)$ is surjective and 
\[
\tx{im}(\pi_{2}(i))\cong \Z_{12}/\tx{ker}(\pi_{2}(i))
=\Z_{12}/\tx{im}(\delta^{k}_{2})=\Z_{12}/k\Z_{12}\cong
\Z_{\tx{gcd}(k,12)}.\qedhere
\]
\end{prf}
\begin{corollary}
If $\cP_{k}$ is a smooth principal $\SU$-bundle over $\bS^{4}$
with Chern number $k$, then
$\pi_{2}(\Gau^{\infty}(\cP))\cong \Z_{\tx{gcd}(12,k)}$, where
$\Gau^{\infty}(\cP)$ denotes the group of smooth gauge transformations
on $\cP$.
\end{corollary}
\begin{prf}
This is the preceding proposition and \cite[Th. III.11]{togg}
\end{prf}
\begin{remark}\label{rem:notationForSurfaces}
We recall that a closed compact orientable surface $\Sigma$ of genus $g$ with
$\partial \Sigma =\emptyset $ can be described as a $\tx{CW}$-complex
by starting with a bouquet
\[
B_{g}=\underbrace{\bS^{1}\vee \dots \vee \bS^{1}}_{2g}
\]
of $2g$ circles. We write $a_{1},b_{1},\dots ,a_{g},b_{g}$ for the
corresponding generators of the fundamental group of $B_{g}$, which is
a free group of $2g$ generators \cite[Th. III.V.14]{bredon}. Then
we consider a continuous map $f:\bS^{1}\to B_{g}$ representing
\[
a_{1}\cdot b_{1}\cdot a_{1}^{-1}\cdot b_{1}^{-1}\cdots
a_{g}\cdot b_{g}\cdot a_{g}^{-1}\cdot b_{g}^{-1}\in \pi_{1}(B_{g}).
\]
Now $\Sigma$ is
homeomorphic to the space obtained by identifying the points on
$\partial \B^{2}\cong \bS^{1}$ with their images in $B_{g}$ under $f$,
i.e. 
\begin{equation}\label{1:eq:descriptionOfSurfaces}
\Sigma\cong B_{g}\cup_{f}\partial \B^{2}
\end{equation}
and we denote by $q_{\Sigma}$ the corresponding quotient map
$q_{\Sigma}\from\B^{2}\to \Sigma$.
\end{remark}
\begin{remark}
Let $\pfb[\Sigma]$ be a continuous principal $K$-bundle over a closed,
compact and orientable surface with, $K$ be connected, and denote by
$q_{\Sigma}\from \B^{2}\to \Sigma$ the quotient map from Remark
\ref{rem:notationForSurfaces}. Then \cite[Cor. VII.6.12]{bredon}
provides a map $\sigma \from \Sigma \to P$ satisfying $\eta \circ
\sigma =q_{\Sigma}$ and since $\left.\cP\right|_{\eta^{-1}(B_{g})}$ is
trivial, we have a continuous map $\gamma\from \partial \B^{2}\to K$
satisfying $\sigma (x)\cdot \gamma (x)=\sigma (y)\cdot \gamma (y)$ if
$x,y\in\partial \B^{2}$ and $f(x)=f(y)$. We may also require
w.l.o.g. that $\gamma(x_{0})=e$ and then $[\gamma]$ may be viewed as a
representative of $\cP $.

Denote by $\sigma'\from B_{g}\to P$ a continuous section. Then $p\sim
p'$ wherever $p=\sigma'(x)\cdot k$ and $p'=\sigma'(y)\cdot k$ for some
$x,y\in B_{g}$ and $k\in K$ defines an equivalence relation on $P$. Then
$P/\sim$ is isomorphic to $P/\gamma$ from Remark
\ref{rem:classificationOfBundlesOverSpheres} (by a similar
construction) and we thus set $P/\gamma :=P/\sim$.
\end{remark}
\begin{proposition}\label{prop:reductionToBundlesOverSpheres}
Let $\cP=\pfb$ be a continuous principal $K$-bundle over a closed,
compact orientable surface, let $K$ be locally contractible and
connected and let $b\in\pi_{1}(K)$ be characteristic for $\cP$ (cf.
Remark
\ref{rem:classificationOfBundlesOverSurfaces}). If $\tx{ev}_{p_{0}}\from
C(P,K)^{K}\to K$ is the evaluation fibration at the base-point of
$P$, then we have a long exact sequence
\begin{multline}\label{eqn:longExactHomotopySequenceForSurfaces}
\begin{CD}\dots 
\to \pi_{n+1}(K)
@>{\delta_{n+1}}>> \pi_{n+1}(K)^{2g}\oplus \pi_{n+2}(K)
\to\pi_{n}(C(P,K)^{K})\end{CD}\\
\begin{CD}\to\pi_{n}(K)
@>{\delta_{n}}>>\pi_{n}(K)^{2g}\oplus \pi_{n+1}(K)
\to \dots
\end{CD}\end{multline}
with connecting homomorphisms $\delta_{n}:\pi_{n}(K)\to
\pi_{n}(K)^{2g}\oplus \pi_{n+1}(K)$ given by $a\mapsto
(0,-\langle a,b \rangle)$, where $\langle \cdot,\cdot \rangle$ denotes the
Samelson product.
\end{proposition}
\begin{prf}
Recall the notation for surfaces from Remark
\ref{rem:notationForSurfaces} and consider the restriction map $r\from
C(P,K)^{K}\to C(\eta^{-1}(B_{g}),K)^{K}$.  Furthermore, \cite[Lem.
IV.4]{togg} provides a continuous map $S\from
C_{*}(\eta^{-1}(B_{g}),K)^{K}\cong C_{*}(B_{g},K)\cong
C_{*}(\bS^{1},K)^{2g}\to C_{*}(P,K)^{K}$ satisfying $r\circ
S=\tx{id}_{C_{*}(\eta^{-1}(B_{g)},K)^{K}}$.  Furthermore,
$\left.\cP\right|_{\eta^{-1}(B_{g})}$ is trivial and we may assume $S$
to be defined on $C(\eta^{-1}(B_{g},K))^{K}$ such that $r\circ
S=\tx{id}_{C_{*}(\eta^{-1}(B_{g)},K)^{K}}$ still holds.

Now take the long exact homotopy sequence for the
fibration $\tx{ev}_{p_{0}}\from C(P,K)^{K}\to K$ \cite[Prop.
IV.8]{togg} and recall the construction of the connecting
homomorphism from Remark
\ref{rem:constructionOfTheConnectingHomomorphism}. If $\alpha\from
\B^{n}\to K$ represents $a\in\pi_{n}(K)$ and $A\from
\B^{n}\to C(P,K)^{K}$ is a lift of $\alpha$, then $A':=A\cdot (S\circ
r\circ A)$ is also a lift of $\alpha$ and $A'(x)(p)=e$ holds
for $x\in \partial \B^{n}$ and $p\in \eta^{-1}(B_{g})$.  Hence
$\left.A'\right|_{\partial \B^{n}}$ factors through a map on
$P/\gamma$ (where $\gamma\from\bS^{1}\to K$ is supposed to represent the
equivalence class of $\cP$) and thus represents $(0,-\langle a,b)
\rangle\in \pi_{n}(K)\oplus \pi_{n+1}(K)$ due to Theorem
\ref{thm:connectingHomomorphismIsTheSamelsonProduct}.
\end{prf}
\begin{remark}\label{rem:rationalHomotopy}
In infinite-dimensional Lie theory one often considers (period-)
homomorphisms $\varphi\from\pi_{n}(G)\to V$ for an
infinite-dimensional Lie Group $G$ and an $\R$-vector space $V$, which
we consider here as a $\Q$-vector space.  If $n\geq 1$, then
$\pi_{n}(G)$ is abelian and this homomorphism factors through the
canonical map $\psi\from\pi_{n}(G)\to \pi_{n}(G)\otimes \Q$, $a\mapsto
a\otimes 1$ and
\[
\wt{\varphi}\from\pi_{n}(G)\otimes \Q\to V,\;\;a\otimes x \mapsto
x\;\varphi (a).
\]
It thus suffices for many interesting questions arising from
infinite-dimensional Lie theory to consider the \textit{rational
homotopy groups} $\pi_{n}^{\Q}(G):=\pi_{n}(G)\otimes \Q$ for
$n\geq 1$.

Furthermore, the functor $\otimes \Q$ in the category of abelian
groups, sending $A$ to $A^{\Q}:=A\otimes \Q$ and $\varphi\from A\to B$
to $\varphi^{\Q}:=\varphi \otimes \id_{\Q}\from A\otimes \Q\to
B\otimes \Q$, preserves exact sequences since $\Q$ is torsion free and
hence flat.
\end{remark}
\begin{lemma}\label{rem:ratinalSamelsonProduct}
If $K$ is a (possibly infinite-dimensioanl) connected Lie group, then the
rational Samelson product
\[
\langle \cdot,\cdot \rangle^{\Q}:\pi_{n}^{\Q}(G)\times \pi_{m}^{\Q}(G)\to 
\pi_{n+m}^{\Q}(G),\;\;
a\otimes x,b\otimes y\mapsto \langle a,b \rangle\otimes xy
\]
vanishes.
\end{lemma}
\begin{prf}
Since each connected Lie group is homeomorphic to a compact group and
a vector space, is has finite-dimensional rational homology and thus
the rational Whitedead product in $BK$ vanishes (cf \cite[Prop. 15.15
f.]{felixHalperinThomas01}).  Since the Whitehead product in $BK$ and
the Samelson product in $K$ correspond to each other via the
connecting homomorphism from the classifying bundle $EK\to BK$
(cf. Remark \ref{rem:alternativeProoOfMainTheorem} and
\cite[Sect. 1]{barrattJamesStein60}), it follows that the rational
Samelson product vanishes either.
\end{prf}
\begin{theorem}
Let $K$ be a connected Lie group and $\cP =\pfb$ be a continuous
principal $K$-bundle over $\bS^{m}$ or a compact orientable
surface $\Sigma$.
\begin{alignat*}{2}
i)\;&\text{If }M=\bS^{m}&\text{, then }&\pi_{n}^{\Q}(\Gau(\cP))\cong
\pi_{n+m}^{\Q}(K)\oplus \pi_{n}^{\Q}(K).\\
ii)\;&\text{If }M=\Sigma&\text{, then }&\pi_{n}^{\Q}(\Gau(\cP))\cong
\pi_{n+2}^{\Q}(K)\oplus \pi_{n+1}^{\Q}(K)^{2g}\oplus \pi_{n}^{\Q}(K).
\end{alignat*}
\end{theorem}
\begin{prf}
With Remark \ref{rem:rationalHomotopy} we obtain exact rational
homotopy sequences from the exact sequences
\eqref{eqn:longExactHomotopySequence} and
\eqref{eqn:longExactHomotopySequenceForSurfaces}. Then the preceding
Lemma implies that the connecting homomorphisms in these sequences
vanish and the long exact sequences decay into short
ones. Furthermore, the short exact sequences split linearly since each of them
involves just vector spaces.
\end{prf}
\begin{remark}
Since the rational homotopy groups of finite-dimensioanl Lie groups
are those of odd-dimensional spheres
\cite[Sect. 15.f]{felixHalperinThomas01}, which are well known
\cite[Ex. 15.d.1]{felixHalperinThomas01} the preceding theorem gives
an explicit formula for $\pi_{n}^{\Q}(\Gau(\cP))$ in the case of
finite-dimensional structure groups. E.g., if $M=\bS^{m}$ and $m$ is
even, then $\pi_{n}^{\Q}(\Gau(\cP))$ vanishes for even $n$.
\end{remark}
\bibliography{mybib}
\vskip\baselineskip \vskip\baselineskip \vskip\baselineskip \large
\noindent
Christoph Wockel\\
Fachbereich Mathematik\\
Technische Universit\"at Darmstadt\\
Schlossgartenstra\ss e 7\\
D-64289 Darmstadt\\
Germany\\[\baselineskip]
\normalsize
\texttt{wockel@mathematik.tu-darmstadt.de}
\end{document}